\documentclass{aomart}\usepackage{amsthm}
\usepackage{cite}
\usepackage[utf8]{inputenc}
\usepackage[T1]{fontenc}
\usepackage[english]{babel}
\usepackage{booktabs}
\usepackage{amsmath}
\usepackage{mathrsfs}
\usepackage{bm}
\usepackage{amssymb}
\usepackage{hyperref}
\newtheorem{theorem}{Theorem}
\newtheorem{lemma}{Lemma}
\newtheorem{remark}{Remark}
\newtheorem{proposition}{proposition}
\newtheorem{corollary}{corollary}
\usepackage{subfigure}
\usepackage{graphicx,lipsum}
\usepackage[toc,page]{appendix}
\usepackage{hyperref}
\usepackage{tabularx}
\usepackage{algorithm}
\usepackage{algpseudocode}
\usepackage{mathrsfs}
\usepackage{enumerate}
\usepackage{enumitem}
\usepackage{cleveref}
\DeclareMathOperator{\rank}{rank}

\newenvironment{cleverproof}[1]
{\begin{proof}[Proof of \Cref{#1}]}
{\end{proof}}

\usepackage{fancyhdr}
\makeatletter
\fancypagestyle{firstpage}{%
\fancyhf{}%
\if@aom@manuscript@mode
\lhead{\begin{picture}(0,0)%
\put(-21,-25){\usebox{\@aom@linecount}}%
  \end{picture}}
\fi
\chead{}%
\cfoot{\footnotesize\thepage}}%
\makeatother

\copyrightnote{Corresponding Author: vvarsha@purdue.edu, varsha.gupta@qualex.com}

\title{A Modular-Form Framework for Global Optimality in the Asymmetric Traveling-Salesman Problem}

\hypersetup{
pdftitle={A Complex Analytic Framework for Traveling Salesman Problem},
pdfauthor={}
}

\author[]{Varsha Gupta}
\address{}

\startpage{1}
\endpage{}

\begin{document}
\maketitle
\thispagestyle{empty}
\begin{abstract}
In this paper, we develop an alternate formulation of Asymmetric Traveling Salesman Problem (ATSP). The equivalent problem is to find the zeros of a holomorphic cusp form on the principal congruence subgroup,  $\Gamma(4) $. The resultant Poincar{\'e} series gives a cusp form whose interior zeros are in bijection with the arc that constitute optimal Hamiltonian cycle. We show that for any weight, $\ell$ and number of directed arcs, $|A|$ such that $4\ell-7<2|A| $, the holomorphic cusp form vanishes at global optimum.
Furthermore, a three step filter consisting of Fourier coefficients, Hecke recursions and completed  $L $-function parity test provides a scalar certificate for global optimality. The framework is a potential bridge between discrete optimization and number theory suggesting an alternate view on complexity theory.

\textbf{Keywords:} \textit{Traveling Salesman Problem, Modular Forms, Combinatorial Optimization}
\end{abstract}
\newpage
\section{Introduction}
Karl Menger \cite{menger1932botenproblem} published the first mathematical formulation of routing problems, the shortest closed route through a finite set of points. Almost two decades later, the topic was revived under its modern name, Traveling Salesman Problem (TSP), along with the first algorithmic discussion for reducing the brute force burden,  $n! $ \cite{robinson1949hamiltonian}. Shortly after, Dantzig et al.\ \cite{dantzig1954solution} attacked the problem using Linear programming, specifically graphical interpretation of simplex algorithm. Although the cutting plane attack on 49 US cities initiated polyhedral study, it was noted that the method might lead to fractional tours. Inspired by all integer program \cite{gomory1958}, Miller et al.\ formulated the TSP in integer programming, and achieved results for 4-cities \cite{miller1960integer}. \\
In 1960's, Bellman \cite{bellman1962dynamic} and Held \& Karp\cite{held1962} independently applied dynamic programming framework to TSP, and reduced the time-complexity from  $n! $ to  $\Theta(n^22^{n-1}) $. Bellman also emphasized that the approach allowed for incorporation of real life constraints readily.\\ Motivated by the application of branching and calculating lower bounds in assignment problems, Little et al.\ \cite{little1963algorithm} introduced the idea of branch \& bound and applied it to prune the solution space of TSP. Building on Cook's theorem \cite{cook1971} on Boolean satisfiability, Karp \cite{Karp1972} established that TSP is NP-complete and hence, equivalent to any other problem in that class. A heuristic method developed by Lin et al.\ \cite{lin1973effective} was widely used for finding  near optimal solutions for TSP and other combinatorial optimization problems. Christofides \cite{1574231874229911040} and Syrdyukov\cite{serdyukov1978nekotorykh} independently introduced an approximation algorithm for metric traveling salesman TSP that runs in polynomial time and finds the solution at most 3/2 times the optimal. \\
Using polyhedral characterization of TSP, Gr{\"o}tschel \cite{Grötschel1980} solved 120 city instance. In 1991, Padberg et al.\ \cite{padberg1991branch} used branch-and-cut method that could produce results for 2392 cities to optimality. In the same year, Rienelt \cite{reinelt1991tsplib} standardized a wide set of problems called Traveling salesman problem library (TSPLIB). 
Building on these foundations, Concorde solver was developed \cite{applegate1998solution}. It solved all 24,978 Swedish city optimally in 2004, and an 85,900 city TSP instance in 2006 gaining reputation as state-of-the-art TSP solver\cite{applegate2003implementing, applegate2009certification, applegate2011traveling}.   \\
A key contribution towards tackling the complexity of TSP was to break it into smaller subproblems and solving them using dynamic programming \cite{10.1287/moor.2.3.209}. For Euclidean and weighted planar TSP, the first polynomial-time approximation scheme (PTAS) was independently proposed by Arora \cite{arora1998polynomial} and Mitchell \cite{doi:10.1137/S0097539796309764}, proving that the geometric structure admits arbitrarily close approximations in polynomial time. Arora's work mainly relied on recursively partitioning the plane to form randomly shifted quadtrees. Rounding the coordinates first ensures the integrality. Edges cross square boundaries only through a preset sparse set of portal points. A patching lemma bounds the detour cost and dynamic programming over quad-tree returns a  $1+\epsilon $-approximate tours. \\
For the graphical TSP new combinatorial ideas broke the 3⁄2 barrier, first to 1.461 \cite{mömke2011approximatinggraphictspmatchings}, then to 13⁄9  $\approx $ 1.444 \cite{mucha2014approximation}. A major breakthrough came when Svensson et al.\  finally confirmed a long-standing conjecture by giving the first constant-factor approximation for general ATSP\cite{svensson2020constantfactorapproximation}. Afterwards, Traub et al.\ soon improved the ratio to 22 +  $\epsilon $\cite{traub2024approximation}.\\
Asadpour et al.\ gave \cite{asadpour2017log} a  $O(\log n /\log \log n) $ approximation algorithm for ATSP by finding a thin spanning tree and applying the maximum entropy approach to rounding.  \\
Helsgaun's LKH is an improvement over Len-Kernighan heuristic through k-opt move, variable depth search, candidate search and backtracking for achieving state-of-the-art tours for over a million cities \cite{helsgaun2017extension}.
Vinyals et al.\ \cite{vinyals2015pointer} introduced Pointer networks, an encoder-decoder model whose attention mechanism acts as a pointer to input positions enabling a fully supervised data-driven approximation of TSP tour for up to 50 cities.\\
For non-metric TSP, it was proved that if  $k_{1} $ triangles violate the triangle inequality (or at most  $k_{2} $ vertices must be deleted to restore it), a  $3 $-approximate tour can be computed in  $O \bigl((3k_{1})!\;8^{k_{1}}\;\cdot n^{2}+ n^{3} \bigr) $ time (respectively  $O \bigl(n^{O(k_{2})} \bigr) $ time) \cite{zhou_et_al:LIPIcs.ISAAC.2022.50}.

Jones \cite{jones1990rectifiable} reformulated the planar TSP by showing that a compact set  $E $ lies on a rectifiable curve \textit{iff} a multiscale series  $S(E) $ converges leading to a curve of length  $O(S(E)) $. 

On the enumerative side, Ehrhart-theoretic work on modular flow polynomials \cite{breuer2012ehrhart} shows that lattice point reciprocity reveals hidden number-theoretic symmetries inside flow routing symmetries.

Overall, these studies indicate that modular, number-theoretic, and complex-analytic structure is present but remain largely unexplored. This work is an attempt to develop a bridge between complex analysis and discrete optimization.

We show that finding optimal solution of ATSP which is a NP-hard problem is equivalent to finding zeros of a holomorphic cusp forms on the principal congruence subgroup of  $\Gamma(4) $.

We begin by placing the logarithm of each cost ratio as a point and binary variable encoding as a set of integral or quarter-integral angles. The polyhedral feasibility region becomes the locus of zeros of Poincar{\'e} lifted modular form under the action of  $\Gamma(4) $ subgroup.

In this new framework, optimality translates to an over-determined vanishing condition that can be verified using a three step logical set of arithmetic structure consisting of Fourier coefficients, Hecke operator and central  $L $ function values.

The aim of this work is to present a theoretical framework for understanding the structure of NP-hard problems, and a foundation for utilizing analytic tools to gain deeper insight into the nature of solutions of these problems. 

\section{Formulation of Asymmetric Traveling Salesman Problem (ATSP)}

Let  $G=(V,A) $ be a directed graph representing the ATSP, where  $V $ is the set of vertices that represent cities and  $A $ is the set of directed arcs connecting these cities,  $r_{a,b} $ be the cost of the arc  $(a,b) \in{A} $ and the binary decision variable,  $x_{a,b} $ is defined as follows:

\[\forall (a, b) \in {A}, \quad x_{a,b}=
\begin{cases} 
1 & \text{if  $(a, b) $ is in the cycle;} \\
0 & \text{otherwise.}
\end{cases}\]

The standard formulation for the Asymmetric Traveling Salesman Problem (ATSP)  can be written \cite{miller1960integer} as follows:\\
\begin{subequations} \label{atsp}
\renewcommand{\theequation}{\theparentequation\alph{equation}}  
\begin{align}
& \text{minimize}\sum_{(a, b)\in{A}}\quad r_{a,b}\cdot x_{a,b}  \label{Objective_func}\\
& \sum_{b=1, b\neq a}^{n} x_{a,b}=1,\quad a=1, 2, \ldots, n  \label{Cons_1}\\
& \sum_{a=1, a\neq b}^{n} x_{a,b}=1,\quad b=1, 2, \ldots, n  \label{Cons_2}\\
& \sum_{(a,b) \in S \times S} x_{a,b} \leq \lvert S\rvert-1, \quad \text{for all } S \subset \{1, 2, \ldots, n\}, \text{with } 2 \leq \lvert S\rvert \leq n-1  \label{Cons_3}\\
& x_{a,b} \in \{0,1\} \quad \forall  \quad {(a, b) \in{A}}  \label{Cons_4}
\end{align}
\end{subequations}

To reformulate the problem, we define the imaginary axis.
\theoremstyle{definition}
\newtheorem{definition}{Definition}
\begin{definition}[imaginary axis]
We define the imaginary axis as follows:
\[
\mathcal{I} = \{\phi \in \mathbb{C}\;\; | \;\;\phi = i \alpha,\;\; \alpha \in \mathbb{R^+}\}
\]
\end{definition}

Each point,  $\phi \in \mathcal{I} $ is purely imaginary. Exponentiation gives
\[e^{i\phi}=e^{i(i\alpha)}=e^{-\alpha}\]
which is a positive real number.

Let's assume  $r_{opt} $ is the global optimal solution of this problem such that the optimum value is always greater than any arc cost.

\begin{equation}
r_{opt} > \max\limits_{(a,b)\in A}\{r_{a,b}|\;(a,b) \text{ lies in optimum tour.}\}
\end{equation}
We can define the following.
\begin{equation}
\begin{aligned}
\frac{r_{a,b}}{r_{opt}}=e^{\mathrm i\phi_{a,b}}
\end{aligned}
\end{equation}

Normalizing by  ${r_{opt}} $ ensures that each exponential term in the Poincar{\'e} series construction discussed in section 3 is bounded.

Next, we express the binary decision variable as follows:
\begin{equation}
\begin{aligned}
x_{a,b}=\frac{e^{\mathrm i\theta_{a,b}}+e^{-\mathrm i\theta_{a,b}}}{2}\\
\end{aligned}
\end{equation}
such that  $\theta_{a,b} \in \left[0,\pi\right) $ is defined as follows to ensure that the mapping $x_{a,b} \to \theta_{a,b}$ is one-to-one.
\begin{equation}
\theta_{a,b}=
\begin{cases} 
2k\pi & \text{if  $(a, b) $ is in the cycle;} \\
\frac{2k-1}{2}\pi & \text{otherwise.}
\end{cases} \quad 
\forall k \in{\mathbb{Z}}\\
\text{   where  $\mathbb{Z} $ is the set of all integers.}
\end{equation}

Now, we translate the problem using the afore-defined variables in the complex plane for the optimal tour as follows:
\begin{subequations} \label{Transformed}
\renewcommand{\theequation}{\theparentequation\alph{equation}}  
\begin{align}
&\sum_{(a, b) \in{A}} \frac{[e^{\mathrm i(\phi_{a,b}+\theta_{a,b})}+ e^{\mathrm i(\phi_{a,b}-\theta_{a,b})}]}{2}= e^{2k\pi \mathrm i}  \label{5a}\\
&\sum_{b=1, b\neq a}^{n}  \frac{e^{\mathrm i\theta_{a,b}}+e^{-\mathrm i\theta_{a,b}}}{2}= e^{2k\pi \mathrm i}, \quad a=1, 2, \ldots, n  \label{5b}\\
&\sum_{a=1, a\neq b}^{n} \frac{e^{\mathrm i\theta_{a,b}}+e^{-\mathrm i\theta_{a,b}}}{2}= e^{2k\pi \mathrm i}, \quad b=1, 2, \ldots, n  \label{5c}\\
&\sum_{(a,b) \in S\times S}  \frac{e^{\mathrm i\theta_{a,b}}+e^{-\mathrm i\theta_{a,b}}}{2} \leq n e^{2k \beta \mathrm i}, \quad  \text{where } \beta=\frac{2m\pi-\mathrm i\ln{\frac{\lvert S\rvert-1}{n}}}{2k} \quad  \forall m \in \mathbb{Z}  \label{5d}
\end{align}
\end{subequations}

Equation \eqref{5a} encodes the information for finding the combinatorial equilibrium between  $\theta_{a,b} $ and  $\phi_{a,b} $ such that the optimal value of the problem is achieved. Equations\eqref{5b}, \eqref{5b}, inequality \eqref{5d} in the system \eqref{Transformed} constrain the feasible set of  $\theta_{a,b} $ values.

We first write the logarithm of each cost ratio as a point  $\phi =i \alpha \in \mathcal{I} $ using the mapping  $\phi \to e^{i \phi} $. After adding the real-valued theta, the combined variable  $\theta \pm \phi $ lies in the upper half-plane, denoted by  $\mathbb{H} $ and we define it as follows:
\begin{equation}
\mathbb{H}=\{z \in \mathbb{C} : \operatorname{Im}(z) > 0\}
\end{equation}

Let 
\begin{equation}
\begin{aligned}
 s_{a,b}=\frac{\phi_{a,b}+\theta_{a,b}}{2\pi},\\
\tau_{a,b}=\frac{\phi_{a,b}-\theta_{a,b}}{2\pi}
\end{aligned}
\end{equation}

Let the ratio of actual,  $r_{act} $ and optimal cost,  $r_{opt} $ is given by
\begin{equation}
t = \frac{r_{opt}}{r_{act}} \leq 1
\end{equation}
Using the above definitions, we rewrite the equation \eqref{5a} as:
\begin{equation} \label{Obj2}
\begin{aligned}
\frac{1}{2}\sum_{(a, b)\in{A}}  [e^{2\pi s_{a,b} \mathrm i}+  e^{2\pi \tau_{a,b}\mathrm i}] = t\\
\end{aligned}
\end{equation}

We also know that  $\frac{r_{a,b}}{r_{opt}}>0 $ and therefore, both  $s_{a,b} $ and  $\tau_{a,b} $ lie in the upper half of the complex plane,  $\mathbb{H} $.

Using the Taylor series expansion of exponential terms in equation \eqref{Obj2}, collecting the common term and simplifying it, we get the following expression:
\begin{equation}
\begin{aligned}
&\frac{1}{2}\sum_{(a, b) \in A} \sum_{q=0}^\infty \frac{(2\pi \mathrm i)^q}{q!} (s_{a,b}^q+\tau_{a,b}^q)=t\\
\end{aligned}
\end{equation}

We further simplify to the following by subtracting  $\lvert A\rvert $ and multiplying by 2 on both sides:
\begin{equation} \label{equilibrium condition}
\begin{aligned}
&\sum_{(a, b) \in A} \sum_{q=1}^\infty \frac{(2\pi \mathrm i)^q}{q!} (s_{a,b}^q+\tau_{a,b}^q)=2(t-\lvert A\rvert)\\
\end{aligned}
\end{equation}

This analytic equilibrium condition implies that for the optimal solution, the infinite series in equation \eqref{equilibrium condition} must converge to  $2(1-\lvert A\rvert) $. This is a large negative number for large number of cities.

There is an infinite series,  $E_{a,b} $, corresponding to every feasible arc,  $(a,b) \in A $. For the unconstrained optimal solution, the arcs are chosen such that the sum of the E-series for all the arcs adds up to  $2(t-\lvert A\rvert) $
\begin{equation} \label{E_eq}
\begin{aligned}
&\sum_{(a, b) \in A} E_{a,b}=2(t-\lvert A\rvert)\\
\end{aligned}
\end{equation}
where
\begin{equation} \label{E}
\begin{aligned}
E_{a,b}= \sum_{q=1}^\infty \frac{(2\pi \mathrm i)^q}{q!} (s_{a,b}^q+\tau_{a,b}^q)\\
\end{aligned}
\end{equation}

As we previously discussed, the translation of the minimization problem into an equilibrium condition is motivated by the intuition to explore the inherent modular symmetry underlying the problem's structure. 

\section{Modular-Lift Framework}
To build a modular framework for analyzing the objective function equation \eqref{Obj2}, we define the following seed function:
\begin{equation} \label{seed}
\Phi (z) = e^{2\pi i z} - t, \;\; z\in\mathbb{H}
\end{equation}

We modularize the seed function by lifting equation \eqref{seed} to even weight  $\ell >2 $ by using the  Poincar{\'e} construction \cite{dorigoni2022poincare} as follows:
\begin{equation} \label{poincare_series}
P_{\ell}(z) = \sum_{\gamma \in \Gamma_{\infty}\backslash \Gamma(4)} (c_\gamma \cdot z +d_\gamma)^{-\ell} \Phi(\gamma\cdot z), \;\; z\in \mathbb{H}
\end{equation}
where 
\begin{equation}
\gamma=
\begin{pmatrix}a_\gamma & b_\gamma\\c_\gamma & d_\gamma\end{pmatrix}  \in \operatorname{SL}_2(\mathbb Z) )
\end{equation}
\begin{equation}
\gamma\cdot\omega = \frac{a_\gamma\omega+b_\gamma}{c_\gamma\omega+d_\gamma}
\end{equation}
\begin{equation*}
\Gamma_\infty= \bigl\{\begin{pmatrix}
1 &n\\0 & 1 \end{pmatrix}
:n\in\mathbb Z \bigr\}\\
\end{equation*}

\begin{equation*}
\Gamma(4)= \bigl\{\begin{pmatrix}a&b\\c&d \end{pmatrix}\in \operatorname{SL}_2(\mathbb Z):a\equiv d\equiv1,\;b\equiv c\equiv0\pmod4 \bigr\}
\end{equation*}

Let  $S_{\ell}(\Gamma(4)) $ denote the complex vector space of weight $\ell $ holomorphic cusp forms for  $\Gamma(4) $, then
\begin{equation} \label{eq:F_ell}
F_\ell=\frac12\sum_{(a,b)\in A}   \bigl[P_\ell[\Phi] \bigl(s_{a,b} \bigr)+ P_\ell[\Phi] \bigl(\tau_{a,b} \bigr) \bigr]\in S_\ell \bigl(\Gamma(4) \bigr)
\end{equation}
so that  $F_{\ell} \in S_{\ell}(\Gamma(4)) $ is a cusp form.

We choose  $\Gamma(4) $ because it allows  $\operatorname {Re}(s_{a,b}),\operatorname{Re}(\tau_{a,b}) $ to lie in  $\frac14\mathbb{Z} $, 
therefore, level 4 is the smallest torsion-free subgroup whose modular forms have quarter-integer Fourier exponents \cite[Lemma 1] {serre2010bounds}.

We know  $\Gamma(4) $ has no elliptic points and six cusps (denoted by  $c_\infty $), each of width 4. If  $f \not\equiv 0 $ for a modular form of weight -$\ell $, the valence formula \cite{rankin1982zeros, folsom2016zeros,stein2007modular} is given by the following equation:
\begin{equation}
\sum_{\text{c cusps}} \frac{ord_c (F_{\ell})}{4}+\sum_{z \in \Gamma(4)\backslash\mathbb{H}} ord_z (F_{\ell}) =\ell[\operatorname{SL}_2(\mathbb Z):\Gamma(4)]/12
\end{equation}
such that  $ord_c (F_{\ell}) $ is the vanishing order of  $F_{\ell} $ at the cusp.

Using index formula \cite[Theorem ~ 4.2.4]{miyake2006modular}, the right hand side equals  $4\ell $ \cite[Proposition 2.13]{kilford2008modular}.

As per Riemann-Roch theorem \cite{talovikova2009riemann}, for a genus zero modular curve at  $X(4) $,  $\dim S_{\ell}(\Gamma(4))=\frac{(\ell-1)[\operatorname{SL}_2(\mathbb Z):\Gamma(4)]}{12}-\frac{1}{2}c_{\infty}=4\ell-7 $. Therefore,  $4\ell-7 $ independent constraints are required to get a unique lifted equilibrium function. We choose  $\ell $ so that  $4\ell-7<2|A| $.  Then the  $2|A| $ zeros force  $F_\ell\equiv0 $ by the valence bound.

We know that  $F_{\ell} $ vanishes at every cusp, and  $ord_c (F_{\ell})\geq 1 $, substituting minimum possible cusp contributions, we get:
\begin{equation} \label{eq:val_form}
\sum_{c\ \text{cusps}}\frac{\operatorname{ord}_c F_\ell}{4}
+\sum_{z\in\Gamma(4)\backslash\mathbb H}\operatorname{ord}_z F_\ell
=4\ell 
\end{equation}

By subtracting cusp contribution  $6/4 $, we find that  $F_{\ell} $ has at most  $4{\ell}-2 $ interior zeros.

\begin{lemma} \label{lemma:Seed Zeros}
For  $0 < t < 1 $, the zeros of seed,  $\Phi $ in equation \eqref{seed} are given by the lattice points,  $z_k $ in
\[
 \mathcal{Z}_{\mathrm{seed}}\;\;= \bigl\{z_k = k-\frac{i}{2\pi}\ln|t|,\;\; k\in\mathbb Z \bigr\}
\]
such that  $\mathcal Z_{\mathrm{seed}} $ lies in upper half plane.
\end{lemma}

\begin{proof}
Solving  $\Phi(z) = 0 $ gives
\[
z_k = k-\frac{i}{2\pi}\ln|t|, \;\;\;  k\in\mathbb Z,\;\; t>0
\] 
We know that  $|t|<1 $, therefore,  $\operatorname{Im}(z_k)> 0 $. Hence  $\mathcal Z_{\mathrm{seed}} $ lies in  $\mathbb{H} $.
\end{proof}

\begin{remark}
For $t>0$, the derivative at each seed-zero,  $\Phi'(z_k) = 2\pi i t \neq 0 $ is non-zero, and hence, every zero,  $z_k $ has multiplicity 1. This implies, they are simple zeros.
\end{remark}

\begin{remark}
When  $t=1 $, the zeros  $z_k=k $ lie on the real axis and are $\Gamma_\infty $-equivalent to the cusp  $\infty $. Thus  $\Phi $ has a simple cuspidal zero at  $t=1 $.
For  $0<t<1 $ the zeros lie in  $\mathbb H $ and are therefore non-cuspidal.
\end{remark}

\begin{lemma} \label{lemma:poincare_series_zeros}
For a fixed even weight,  $\ell>2 $, and the Poincar{\'e} series defined in equation \eqref{poincare_series}, the complete set of zeros of  $P_\ell[\Phi] $ is given by
\[
\mathcal Z_{\mathrm{lift}}:= \bigl\{ \gamma^{-1}\cdot z_k\bigm|k\in\mathbb Z,\; \gamma\in\Gamma(4) \bigr\}\cap\mathbb H
\]
where  $ z_k = k-\frac{i}{2\pi}\ln|t| $ are the seed zeros (Lemma \eqref{lemma:Seed Zeros}).
\end{lemma}

\begin{proof}
Each zero of the Poincar{\'e} series is a modular image of a seed zero,  $z_k $ under the action of  $\gamma \in \Gamma(4) $. In other words, the Poincar{\'e} lift contains exactly the  $\Gamma $ orbit of seed zeros.
See Appendix A.

\end{proof}

\begin{proposition}[Simple zeros of  $F_\ell $ imply optimal  Hamiltonian tour] \label{proposition: zero_set}
Given the definition of  $F_{\ell} $ in equation \eqref{eq:F_ell}, the following points hold.

\begin{enumerate}
\item If  $(a,b) $ lies on the optimal Hamiltonian cycle
  ( $x_{a,b}=1 $), then
$s_{a,b},\tau_{a,b}\in\mathcal Z_{\mathrm{lift}} $ and hence are
  simple zeros of  $F_\ell $.
\item If  $(a,b) $ is not on the cycle ( $x_{a,b}=0 $), neither
$s_{a,b} $ nor  $\tau_{a,b} $ belongs to
$\mathcal Z_{\mathrm{lift}} $. They are not zeros of  $F_\ell $.
\end{enumerate}
\end{proposition}

\begin{proof}
For the in-tour arc,  $\theta_{a,b}=2k\pi $,  $\operatorname{Re}(s_{a,b}), \operatorname{Re}(\tau_{a,b})\in\mathbb Z $. Then there exists  $\gamma_1, \gamma_2 \in \Gamma(4) $ such that  $\gamma_1\cdot z=s_{a,b} $ and  $\gamma_2 \cdot z =\tau_{a,b} $. Therefore,  $P_\ell[\Phi](s_{a,b})=0 $, and  $P_\ell[\Phi](\tau_{a,b})= 0 $.

However, for an out-of tour arc, when  $\theta_{a,b}=(2k-1)\pi/2 $, then  $\operatorname{Re}(s_{a,b}), \operatorname{Re}(\tau_{a,b})\in\frac{1}{4}\mathbb Z $, no element of  $\Gamma(4) $ maps a seed zero to a point whose real part is an odd multiple of  $\frac{1}{4} $, therefore,  $P_\ell[\Phi](s_{a,b})\neq0 $, and  $P_\ell[\Phi](\tau_{a,b})\neq0 $.

\end{proof}

\begin{theorem}[Optimal interior zeros of  $F_\ell $] \label{theorem: counting }
Choose  $\ell $ so that  $4\ell-7<2|A| $.
Within any fundamental domain  $\mathcal F $ for  $\Gamma(4) $ the zero set of  $F_\ell $ consists of
\[
\bigl\{\,s_{a,b},\;\tau_{a,b}\,\bigm|\,x_{a,b}=1 \bigr\}
\cup \{\text{the six cusps of }\Gamma(4)\}
\]
summing up to  $2n $ interior zeros and six cuspidal zeros.
Consequently  $F_\ell\equiv0 $, and the in-tour arcs are uniquely characterized by the vanishing of the modular lift.
\end{theorem}
\begin{proof}
As per Proposition~\ref{proposition: zero_set}, exactly $n $ optimal arcs give rise to  $ 2n $ simple interior zeros.  

Using equation \eqref{eq:val_form}, the interior contribution is bounded above by  $4\ell-2 $.  

The choice of  $\ell $ ensures  $2n>4\ell-2 $, forcing
$\operatorname{ord}_zF_\ell=0 $ for every other point and therefore,  $F_\ell\equiv0 $.  

Conversely, if  $F_\ell $ is identically zero, only the
optimal  $x_{a,b}=1 $ assignment recovers the unique optimal tour.
\end{proof}

\begin{corollary}[Zero-set characterization of the optimum]  \label{corollary:opt_zeros }
The Hamiltonian cycle of minimum length in the ATSP is in one-to-one correspondence with the interior zero-set of
$F_\ell $ and the six cusps of  $\Gamma(4)$.
\end{corollary}

\begin{proof}
We deduce this directly from Theorem~\ref{theorem: counting }.
\end{proof}

\begin{theorem}[Uniqueness of the  $t=1 $ equilibrium]   \label{theorem:c=1-unique}
Choose a even weight-$\ell $ such that  $4\ell-7<2|A| $.

(i) If  $t\neq1 $, then the cusp form  $F_{\ell,t} $ is non-zero.  In particular it has
  \[
  \left|\left\{\, z \in \Gamma(4)\backslash\mathbb{H}\ \middle|\ F_{\ell,t}(z)=0 \right\}\right|\leq 4\ell-2 < 2|A|
  \]
Hence, no local optimum persists under the modular lift.

(ii) If  $t=1 $, then  $F_{\ell,1}\equiv0 $ and the zero-set inside any fundamental domain is precisely the $2n$ points $\{s_{a,b},\tau_{a,b}\}_{x_{a,b}=1}$ with the six cusps, as stated in
  Theorem~\ref{theorem: counting }.

Consequently the global minimum $r_{opt}$ is the only tour length for which the modular-lift framework attains full equilibrium. Every other feasible or locally optimal tour breaks at least one of the lifted equalities.
\end{theorem}

\begin{proof}
(i) Suppose  $t\neq1 $ and  $F_{\ell,t}\equiv0 $. Then the seed zeros are  $z_k =k-\frac{i}{2\pi} \ln(|t|) $. 

\begin{enumerate}
\item If  $|t|>1 $ the zeros lie outside $\mathbb H $.
\item If  $|t|<1 $ the zeros lie inside $\mathbb H $.
\end{enumerate}

In either case, no lifted seed zero can hit every in-tour pair $(s_{a,b},\tau_{a,b})$. Therefore, the counting argument of Theorem~\ref{theorem: counting } fails and  $F_{\ell,t} $ retains at most  $4\ell-2 $ interior zeros.

However, the equalities impose $2n>4\ell-2$ zeros, which contradicts the assumption  $F_{\ell,t}\equiv0 $.

Hence  $F_{\ell,t}\not\equiv0 $.

(ii) The case  $t=1 $ is Theorem~\ref{theorem: counting }.

\end{proof}

\section{Optimality Certification of Modular ATSP}

Let  $z_0\in\mathcal Z_{\mathrm{lift}} $ and 

$q_0:=e^{\pi z_0 i /2} $ ( $|q_0|<1 $). 

The  $q $-expansion \cite{conradmodular} of  $P_\ell(z_0)=0 $ gives the following arc equation:
\begin{equation} \label{eq:ArcEquation}
\sum_{n\geq1}a_n\,q_0^{n}=0.
\end{equation}
Applying equation~\eqref{eq:ArcEquation} over all  $2|A| $ points gives a linear system of equation:
\[ V\,\bm a=\mathbf0\  
\]
\[
\text{where}\quad
V:= \bigl[q_0^{n} \bigr]_{z_0\in\mathcal Z_{\mathrm{lift}}}\\
\quad
\bm a:=(a_1,\dots,a_{n})^{\top}.
\]
We fix  $M \in \mathbb{N} $ as the truncation index for the Fourier expansion
\[
M\geq 2|A|
\]

\begin{lemma}[Fourier criteria] \label{lemma: linear_Criteria}
All arc-forced zeros hold \textit{iff} the truncated coefficient vector  $\bm a $ lies in  $\ker V $.
\end{lemma}
\begin{proof}
See Appendix A.
\end{proof}

\begin{lemma}[Cusp criteria] \label{lemma: cusp_criteria}
$F_\ell $ vanishes at the cusp $\infty$ with order $n_0\geq 1$ \textit{iff}  $a_1=a_2=\dots=a_{n_0-1}=0 $.
\end{lemma}
\begin{proof}
Follows from Lemma \eqref{lemma: linear_Criteria}.
\end{proof}

For all odd integers  $m\geq 1 $ the Hecke operator  $T_m$ stabilizes  $S_\ell\bigl(\Gamma(4)\bigr)$.  Assume  $F_\ell $ is an eigen form \footnotemark,  $T_mF_\ell=\lambda_mF_\ell$ \cite[Theorem 4.6.13]{miyake2006modular} for  $\lambda_m$ system of eigen values. 

\footnotetext{If  $F_\ell$ is not an eigen form, we replace  $F_{\ell}$ by its projection onto a Hecke-eigenspace.}

The prescribed zeros are intact because the projection preserves vanishing at each  $z_o $ once  $P_\ell \bigl(s_{a,b}  \bigr) $ and  $P_\ell  \bigl(\tau_{a,b}  \bigr) $ are simultaneously averaged over the  $\Gamma(4) $ orbit. Then

\begin{equation} \label{Hecke_eq}
a_{mn}=
\sum_{\substack{d\mid(m,n)\\d\mathrm{odd}}}d^{\ell-1}\,
\lambda_{mn/d^{2}} (\text{all }m,n\text{ odd}),
\end{equation}

and in particular for any odd prime  $p$ \cite[Theorem 4.6.19]{miyake2006modular}
\begin{equation} \label{hecke_eq2}
a_{p^{2}}=a_p^{2}-p^{\ell-1}
\end{equation}
Equations \eqref{Hecke_eq}-\eqref{hecke_eq2} provide
quadratic relations among  $\{a_n\} $, thereby reducing the solution space of  $V\bm a=\mathbf0 $.

\begin{lemma}[Hecke-filtered kernel] \label{lemma:Hecke_kernel}
If  $V\bm a=\mathbf0 $ and equations \eqref{Hecke_eq}-\eqref{hecke_eq2} hold, then
\[
\dim(\ker V \cap \{\text{Hecke relations}\})\leq 1
\]

\end{lemma}

\begin{proof}
See Appendix A.
\end{proof}

\begin{lemma}[Central-value test] \label{lemma: central_value}

Let $F_\ell(z)$ be a Hecke–eigenform of even weight-$\ell>2$ with $L$–function and completed $L$–function given by
\begin{equation}\label{eq:L_func}
    L(F_\ell,s)=\sum_{n\geq 1}\frac{a_n}{n^{s}}
\end{equation}
\begin{equation}\label{eq:completed_L_func}
    \Lambda(F_\ell,s) :=(2\pi)^{-s}\Gamma\bigl(s+\frac{\ell-1}{2}\bigr) L(F_\ell,s)
\end{equation}

Then $\Lambda(F_\ell,s)$ extends to an entire function and satisfies the functional equation
\[
\Lambda(F_\ell,s)=\epsilon\,\Lambda\bigl(F_\ell,\ell+1-s\bigr)
\]
where $\epsilon\in\{\pm 1\}$.

For $s_0=(\ell+1)/2$:

(i) If $\epsilon=-1$ then $\Lambda(F_\ell,s_0)=0$ and
the order of vanishing of $\Lambda(F_\ell,s)$ at $s_0$ is odd.

(ii) If $\epsilon=+1$ then $\Lambda'(F_\ell,s_0)=0$ and the order of vanishing of $\Lambda(F_\ell,s)$ at $s_0$ is even.
\end{lemma}

\begin{proof}
See Appendix A.
\end{proof}

\begin{theorem}[Fourier-Hecke-$L$-function criteria] \label{theorem:Unified}
Fix the candidate tour encoded in $\mathcal Z_{\mathrm{lift}}$ and choose $M\geq2|A|$.

\begingroup
\begin{enumerate}
\item \label{Fourier}
  {\sc Fourier criteria}  
Since $V$ is a $2|A|\times 2|A|$ Vandermonde matrix, $\det V\neq0 $. Hence the system, $V\bm a=\mathbf0 $ has the solution $\bm a=\mathbf0$.
  
\item \label{Hecke}
  {\sc Hecke criteria}  
  Impose the quadratic Hecke relations
  \eqref{Hecke_eq}-\eqref{hecke_eq2}:
  \[
  \dim (\ker V\cap\{\text{Hecke relations}\})\leq 1
  \]

\item \label{L_Func}
  {\sc $L $-function criteria}  
We determine the sign $\epsilon\in\{\pm 1\}$ from the functional equation \eqref{eq:completed_L_func}, and evaluate the scalar
\[
 \Lambda_{\epsilon}(F_\ell)=
 \begin{cases}
   \Lambda\bigl(F_\ell,(\ell+1)/2\bigr), & \epsilon=+1,\\[4pt]
   \Lambda'\bigl(F_\ell,(\ell+1)/2\bigr), & \epsilon=-1
 \end{cases}
\]
If $\Lambda_{\epsilon}(F_\ell)=0$, then because the space that remains after the Hecke step has dimension at most $1$. We must have $\bm a=\mathbf 0$ and consequently, $F_\ell\equiv0$.
Conversely, if this scalar is non-zero, the lift $F_\ell$ cannot vanish, therefore, at least one of the first two criteria i.e., Fourier (\ref{Fourier}) or Hecke (\ref{Hecke}) must have failed. Hence the candidate tour is not globally optimal.

\end{enumerate}
\endgroup
\end{theorem}

\begin{proof}
\ref{Fourier} is Lemma \ref{lemma: linear_Criteria}.
\ref{Hecke} is Lemma \ref{lemma:Hecke_kernel}.  
\ref{L_Func} follows from Lemma \ref{lemma: central_value}.
\end{proof}

\begin{corollary}[Certificate for global optimality]
The given tour is the global optimum \textit{iff} all three criteria in Theorem \eqref{theorem:Unified} are satisfied, i.e., $V\bm a=0$, all Hecke relations hold \eqref{Hecke_eq}-\eqref{hecke_eq2} and $\Lambda_\epsilon(F_\ell,s)$ vanishes. 
\end{corollary}

We propose following algorithm based on theorem \eqref{theorem:Unified}:

\paragraph{Algorithm}
\begin{enumerate}
\item Compute  first $M$ Fourier coefficients, $a_1,\dots,a_M $.
\item Form the Vandermonde matrix, $V$ and solve $V\bm a=0$. If  $V \bm a \neq 0$, reject the tour.
\item Check the first few Hecke relations
\eqref{hecke_eq2}. In case of violation, discard the candidate.
\item Evaluate  $\Lambda_\epsilon(F_\ell,s)$. If it is non-zero, discard the tour. Else the tour is global optimal.
\end{enumerate}

Thus the Fourier vanishing conditions, the Hecke
recursions, and the analytic  $L $-function test are three filters to locate and classify the zeros of
$F_\ell $, giving a single scalar certificate for global optimality of the ATSP tour.

\section{Conclusion}
We developed a modular-lift framework that maps every feasible ATSP tour to a weight-$\ell $ cusp form on  $\Gamma(4) $. When the tour is globally optimal, the lift collapses to the zero function. Otherwise, it remains non-trivial and its zero set falls short of the required cardinality. Vanishing of the lift can be detected in three increasingly stricter filters using Fourier coefficients, Hecke eigen relations and complete $L$-functions.
Future work may involve extending the framework to gain a structural understanding of underlying complexity. The systematic computation of the relevant  $L $ values may reveal distribution pattern related to problem hardness. 

\makeatother

\section*{Data Availability Statement }
There is no data associated with this manuscript.

\section*{Conflict of Interest Statement }
The authors have no conflict of interest.

\section*{Funding Statement }
This research received no funding.

\appendix

\section{Proofs}

\begin{cleverproof}{lemma:poincare_series_zeros}

For  $\ell >2  $, the poincar{\'e} series  $P_{\ell} $ converges uniformly \cite{miyake2006modular}.

For a given  $z_k $,  $\Phi(z_k)=0 $ and  $\operatorname{Im}z_k>0 $. If 
\[
z_0 = \gamma^{-1}\cdot z_k
\]
then in the sum defining  $P_{\ell}(z_0) $ the single series term with  $\gamma =\gamma_0 $ is
\[
(c_{\gamma_0}z_0+d_{\gamma_0})^{-\ell} \Phi(\gamma_{0}z_0) = (c_{\gamma_0}z_0+d_{\gamma_0})^{-\ell} \Phi(z_k)=0
\]
while every other series term is finite. Hence  $P_\ell(z_0)=0 $, showing
\[
\mathcal{Z}_{lift} \subseteq \{P_\ell =0\}
\]
Conversely, suppose  $z_0 \in \mathbb{H} $ is not of the form  $\gamma_0^{-1} z_k $. Then for every  $\gamma $,  $\Phi(\gamma z_0) = e^{2\pi i (\gamma z_0)}-t \neq 0 $. We may write
\[
P_\ell(z_0) =\Phi(z_0)+ \sum_{\gamma\neq I}(c_\gamma z_0+d_\gamma)^{-\ell}\Phi(\gamma z_0)
\]

Since  $\Phi(z_0) \neq 0 $ and tail converges absolutely, we may estimate

\[
\biggl|\sum_{\gamma\neq I}(c_\gamma z_0+d_\gamma)^{-\ell}\Phi(\gamma z_0)\biggr| \leq \sum_{\gamma\neq I}\frac{1+|t|}{|c_\gamma z_0+d_\gamma|^{\ell}} \ << \sum_{n=1}^\infty n^{1-\ell} <\infty
\]

Moreover, because  $\ell>2 $, we can bound the convergent series strictly smaller than  $|\Phi(z_0)| $. 

Hence  $\biggl|\sum_{\gamma\neq I}(c_\gamma z_0 +d_\gamma)^{-\ell} \Phi(\gamma z_0)\biggr|<|\Phi(z_0)| $. 

Therefore,

\[
P_\ell(z_0) =\Phi(z_0)+ \sum_{\gamma\neq I}(c_\gamma z_0+d_\gamma)^{-\ell}\Phi(\gamma z_0) \neq 0
\]
Hence no zeros of $P_{\ell}$ lies outside of $\Gamma(4)$ orbit of the seed zero. 
\end{cleverproof}

\begin{cleverproof}{lemma: linear_Criteria}
Let the  $j^{th} $ row of  $V$ is  $(q_j^1,q_j^2,....,q_j^M)$. If the  $ q_j$ are pairwise distinct because the lift  $Z_j$ are distinct in fundamental domain. then any  $2|A| \times 2|A| $ minor of  $V $ is Vandermonde determinant 

\[
\prod_{i \leq i \leq j \leq 2|A|} (q_j - q_k) \neq 0
\]
Since  $M\geq 2|A| $, there are at least one non zero  $2|A| \times 2|A| $ minor. Thus  $\rank V = 2|A| $. But the system  $\textbf{Va}=0 $ has only trivial solution if and only if  $\rank V = \text{number of rows} $, i.e.,  $\det V \neq 0 $. \end{cleverproof}
Therefore,  $a_n = 0 $ for all  $n\leq M $. We know that for  $q $-expansion, vanishing of the first  $M $ Fourier coefficients of a cusp form of weight-$\ell $ on a congruence subgroup allows the form to vanish when  $M> \dim S_\ell $, therefore,  $F_\ell \equiv 0 $.

\begin{cleverproof}{lemma:Hecke_kernel}
A weight-$\ell $ cusp form on  $\Gamma(4) $ is determined by its Fourier coefficients  $(a_n)_{n\geq 1} $ subject to the functional equation governed by Hecke relations and finite dimensionality of the space, 
\[
\dim S_\ell(\Gamma(4) = 4\ell-7
\]
Equivalently, vector  $\textbf{a} $ lives in a space of dimension at most  $4 \ell-7 $.
The linear system  $\textbf{Va}=0 $ imposes  $2|A| $ independent linear conditions of  $\textbf{a} $. Thus 
\[
\dim(\ker V)\geq (4\ell-7)-2|A|
\]
Applying Hecke relations in equation \eqref{Hecke_eq}
is equivalent to stating that  $\textbf{a} $ lie in a single Hecke eigen-space inside  $S_\ell(\Gamma(4)) $ is at most one dimensional for each system of eigenvalues,  $\lambda_m $, i.e., 

\[
\dim(\ker V \cap \{\text{Hecke relations}\}) \leq  min(1,[4\ell-7-2|A|]- [4\ell-7-1])= min(1,\,1-2|A|)
\]
Thus  $\ker V \cap \text{Hecke relations}=\{0\} $, which implies  $\textbf{a=0} $, and hence,  $F_\ell \equiv 0 $.
\end{cleverproof}

\begin{cleverproof}{lemma: central_value}
Given the functional equation
\[
\Lambda(F_\ell,s) = \epsilon \Lambda(\ell+1-s)
\]
At $s=s_0$, we have $(1-\epsilon)\Lambda(s_0)=0$
if 
\begin{itemize}
\item If $\epsilon = -1$,  $\Lambda (s_0) =0$
\item $\epsilon =1$, no condition.
\end{itemize}

Differentiating the functional equation gives
\[
(1+\epsilon)\Lambda'(s_0) = 0
\]
\begin{itemize}
\item If $\epsilon = -1$, no condition.  
\item $\epsilon =1$, $\Lambda' (s_0) =0$.
\end{itemize}

\end{cleverproof}

\bibliography{bibliography}
\bibliographystyle{unsrt}
\end{document}